# A branch and bound technique for finding the minimal solutions of the linear optimization problems subjected to Lukasiewicz

Amin Ghodousian
Faculty of Engineering Science
University of Tehran
Tehran, Iran
a.ghodousian@ut.ac.ir

Zahra Boreiri
Department of Engineering Science
University of Tehran
Tehran, Iran
zahra.boreiri@ut.ac.ir

Faeze Habibi
Department of Engineering Science
University of Tehran
Tehran, Iran
faeze.habibi@ut.ac.ir

*Abstract*— In this paper, an optimization model with a linear objective function subjected to a system of fuzzy relation equations (FRE) is studied where the feasible region is defined by the Lukasiewicz t-norm. Since the finding of all minimal solutions is an NP-hard problem, designing an efficient solution procedure for solving such problems is not a trivial job. Firstly, the feasible domain is characterized and then the problem is solved with a modified branch-and-bound solution technique based on a new solution set that includes the minimal solutions. After presenting our solution procedure, a concrete example is included for illustration purpose.

## I. INTRODUCTION

In this paper, we study the following linear problem in which the constraints are formed as fuzzy relational equations defined by Lukasiewicz t-norm:

$$\begin{aligned} \min \quad & cx \\ & A\varphi x = b \\ & x \in [0,1]^n \end{aligned} \quad (1)$$

Where $I = \{1,2,\ldots,m\}$, $J = \{1,2,\ldots,n\}$, $A$ is a fuzzy matrix such that $0 \leq a_{ij} \leq 1$ ($\forall i \in I$ and $\forall j \in J$), $b = (b_i)_{m \times 1}$ is a fuzzy vector such that $0 \leq b_i \leq 1$ ($\forall i \in I$), and "$\varphi$" is the Lukasiewicz t-norm defined by $\varphi(x,y) = T_L(x,y) = \max\{x + y - 1, 0\}$. If $a_i$ is the $i$'th row of matrix $A$, then problem (1) can be expressed as follows:

$$\begin{aligned} \min \quad & cx \\ & \varphi(a_i, x) = b_i \,, \, i \in I \\ & x \in [0,1]^n \end{aligned}$$

where the constraints mean:

$$\varphi(a_i, x) = \max_{j \in J}\{\varphi(a_{ij}, x_j)\} = \max_{j \in J}\{T_L(a_{ij}, x_j)\}$$

$$\varphi(a_i, x) = \max_{j \in J}\{\varphi(a_{ij}, x_j)\} = \max_{j \in J}\{T_L(a_{ij}, x_j)\}$$

The theory of fuzzy relational equations (FRE) was firstly proposed by Sanchez and applied in problems of the medical diagnosis [13]. Nowadays, it is well known that many issues associated with a body knowledge can be treated as FRE problems [12]. The solution set of FRE is often a non-convex set that is completely determined by one maximum solution and a finite number of minimal solutions [4]. The other bottleneck is concerned with detecting the minimal solutions that is an NP-hard problem [5,7,8,10]. The problem of optimization subject to FRE and FRI is one of the most interesting and on-going research topics among the problems related to FRE and FRI theory [1,2,4–8,11,14]. Recently, many interesting generalizations of the linear programming subject to a system of fuzzy relations have been introduced and developed based on composite operations used in FRE, fuzzy relations used in the definition of the constraints, some developments on the objective function of the problems and other ideas [3,7-9,11].

In this paper, after carefully study the solution set of system (1), we show that the problem can be converted to a 0-1 integer programming problem and solved by a branch-and-bound method.

The rest of the paper is arranged as follows. In Section 2, a special characterization of the feasible domain of problem (1) is derived. Based on the theoretical aspects of the problem, a new solution set is obtained that includes all the minimal solutions. In Section 3, we study the effect of the cost vector c with the special characterization. A modified branch-and-bound method is presented and a step-by-step algorithm for solving problem (1) is given in Section 4. Some examples are included in Section 5 to illustrate how the algorithm works. Some concluding remarks are made in Section 6.

## II. CHARACTERIZATION OF FEASIBLE REGION

Let $S(A,b)$ denote the feasible solutions set of problem (1), that is,

$$S(A,b) = \left\{ x \in [0,1]^n : \max_{j=1}^{n}\{T_L(a_{ij}, x_j)\} = b_i, \forall i \in I \right\}$$

Also, for each $i \in I$, define $J_i = \{j \in J : a_{ij} \geq b_i\}$. According to [8], when $S(A,b) \neq \emptyset$, it can be completely determined by one maximum solution and a finite number of minimum solutions. The maximum solution can be obtained by $\bar{X} = \min_{i \in I}\{\hat{x}_i\}$ where $\hat{x}_i = [(\hat{x}_i)_1, \ldots, (\hat{x}_i)_n]$ $(i = 1, \ldots, m)$ is defined as follows

$$(\hat{x}_i)_j = \begin{cases} b_i + 1 - a_{ij} & j \in J_i \\ 1 & j \notin J_i \end{cases}, \forall j \in J \quad (2)$$

Moreover, if we denote the set of all minimum solutions by $\underline{S}(A,b)$, then

$$S(A,b) = \bigcup_{\underline{X} \in \underline{S}(A,b)} \{x \in [0,1]^n : \underline{X} \leq x \leq \bar{X}\} \quad (3)$$

**Definition 1.** Let $i \in I$. For each $j \in J_i$, we define $\check{x}_i(j) = [\check{x}_i(j)_1, \ldots, \check{x}_i(j)_n]$ such that

$$\check{x}_i(j)_k = \begin{cases} b_i + 1 - a_{ij} & b_i \neq 0 \text{ and } k = j \\ 0 & otherwise \end{cases}, \forall k \in J$$

**Definition 2.** Let $e: I \to \bigcup_{i \in I} J_i$ so that $e(i) = j \in J_i$, $\forall i \in I$, and let $E$ be the set of all vectors $e$. For the sake of convenience, we represent each $e \in E$ as an $m$–dimensional vector $e = [j_1, \ldots, j_m]$ in which $j_k = e(k)$.

**Definition 3.** Let $e = [j_1, \ldots, j_m] \in E$. We define $\underline{X}(e) = [\underline{X}(e)_1, \ldots, \underline{X}(e)_n]$ where $\underline{X}(e)_j = \max_{i \in I}\{\check{x}_i(e(i))_j\}$, $\forall j \in J$.

The following figure shows the network generated by vectors $e = [j_1, \ldots, j_m]$:

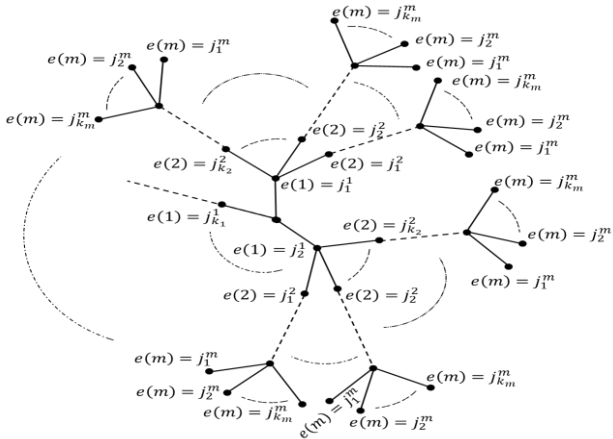

Based on the definitions 2 and 3, and according to [8], we have:

$$\underline{S}(A,b) \subseteq \{\underline{X}(e) : e \in E\} \quad (4)$$

**Definition 4.** For each $i \in I$, define $\bar{J}_i = \{j \in J_i : T_L(a_{ij}, \bar{X}_j) = b_i\}$. Also, let $\bar{E}$ be the set of all vectors $e: I \to \bigcup_{i \in I} \bar{J}_i$ so that $e(i) = j \in \bar{J}_i$, $\forall i \in I$.

**Lemma 1.** Let $i \in I$, $j \in J_i$. If $x_j < b_i + 1 - a_{ij}$, and $b_i \neq 0$, then $T_L(a_{ij}, x_j) < b_i$.

**Proof.** The proof is easily obtained from the definition of Lukasiewicz t-norm. □

**Theorem 1.** Suppose that $e \in \bar{E}$. Then, $\underline{X}(e) \leq \bar{X}$.

**Proof.** According to Definitions 1, 2 and 3, either $\underline{X}(e)_j = 0$ or $\underline{X}(e)_j = b_{i_0} + 1 - a_{i_0 j}$ for some $i_0 \in I$ such that $e(i_0) = j$. In the former case, we have obviously $\underline{X}(e)_j = 0 \leq \bar{X}_j$. In the latter case, $\underline{X}(e)_j = b_{i_0} + 1 - a_{i_0 j}$ and $e(i_0) = j \in \bar{J}_i$. Therefore, by letting $I_j = \{i \in I : a_{ij} \geq b_i\}$, we have $i_0 \in I_j$ and then $\underline{X}(e)_j \geq \min_{i \in I_j}\{b_i + 1 - a_{ij}\} = \bar{X}_j$. To complete the proof, it is sufficient to show that the inequality never holds. By contradiction, suppose that $\underline{X}(e)_j = b_{i_0} + 1 - a_{i_0 j} > \bar{X}_j$. Hence, by Lemma 1 we have $T_L(a_{i_0 j}, \bar{X}_j) < b_{i_0}$ which contradicts $e(i_0) \in \bar{J}_i$. □

**Corollary 1.** $\{\underline{X}(e) : e \in \bar{E}\} \subseteq \{\underline{X}(e) : e \in E\}$.

**Proof.** From Definition 4, we have $\bar{J}_i \subseteq J_i$. Therefore, $\bar{E} \subseteq E$ that implies $\{\underline{X}(e) : e \in \bar{E}\} \subseteq \{\underline{X}(e) : e \in E\}$. □

**Corollary 2.** Suppose that $e \in E - \bar{E}$. Then, $\underline{X}(e)$ is an infeasible point.

**Proof.** If $e \in E - \bar{E}$, then from Theorem 1, the inequality $\underline{X}(e) \leq \bar{X}$ is not satisfied. Now, the result follows from Relation (3). □

**Corollary 3.** $\underline{S}(A,b) \subseteq \{\underline{X}(e) : e \in \bar{E}\} \subseteq \{\underline{X}(e) : e \in E\}$ Suppose that $e \in E - \bar{E}$. Then, $\underline{X}(e)$ is an infeasible point.

**Proof.** The proof is resulted from Corollaries 1 and 2. □

**Corollary 4.** $S(A,b) = \bigcup_{e \in \bar{E}} [\underline{X}(e), \bar{X}]$.

**Proof.** The proof is resulted from Relation (3) and Corollary 3. □

## III. Linear Programming and Effect of Cost Vector

According to the well-known schemes used for optimization of linear problems such as (1) [5,10], problem (1) is converted to the following two sub-problems:

$$\min \quad Z_1 = \sum_{j=1}^{n} c_j^+ x_j$$
$$\varphi(a_i, x) = b_i, \quad i \in I$$
$$x \in [0,1]^n$$

And,

$$\min \quad Z_2 = \sum_{j=1}^{n} c_j^- x_j$$
$$\varphi(a_i, x) = b_i, \quad i \in I$$
$$x \in [0,1]^n$$

where $c_j^+ = \max\{c_j, 0\}$ and $c_j^- = \min\{c_j, 0\}$ for $j = 1, 2, \ldots, n$. It is easy to prove that $\bar{X}$ is the optimal solution of $Z_2$, and the optimal solution of $Z_1$ is $\underline{X}(e')$ for some $e' \in \bar{E}$.

**Theorem 2.** Suppose that $S(A, b) \neq \emptyset$ and $\bar{X}$ and $\underline{X}(e^*)$ are the optimal solutions for $Z_2$ and $Z_1$, respectively. Then $c^T x^*$ is an optimal solution of problem (1), where $x^* = [x_1^*, x_2^*, \ldots, x_n^*]$ is defined as follows:

$$x_j^* = \begin{cases} \bar{X}_j & c_j \leq 0 \\ \underline{X}(e^*)_j & c_j > 0 \end{cases} \quad (5)$$

for $j = 1, 2, \ldots, n$.

**Proof.** Let $x \in S(A, b)$. Then, from Corollary 4 we have $x \in \bigcup_{e \in \bar{E}} [\underline{X}(e), \bar{X}]$. Therefore, for each $j \in J$ such that $c_j > 0$, inequality $x_j^* \leq x_j$ implies $c_j^+ x_j^* \leq c_j^+ x_j$. In addition, for each $j \in J$ such that $c_j \leq 0$, inequality $x_j^* \geq x_j$ implies $c_j^- x_j^* \leq c_j^- x_j$. Hence, $\sum_{j=1}^{n} c_j x_j^* \leq \sum_{j=1}^{n} c_j x_j$. □

As mentioned earlier, generating the maximum solution $\bar{X}$ is not a problem. If we know how to find a minimum solution for $Z_1$, then problem (1) can be solved via Theorem 2. In the next section, we provide a branch and bound method for solving $Z_1$.

## IV. Modified Branch and Bound Method

A branch-and-bound method implicitly enumerates all possible solutions to a programming problem. For our application, in the beginning, we choose one constraint to branch the original problem into several sub-problems. Each sub-problem is represented by one node. Then branching at each node is done by adding one additional constraint. New sub-problems are created and represented by new nodes. Note that the more constraints added to a sub-problem, the smaller feasible domain it has and, consequently, the larger optimal objective value Z it achieves. Therefore, solving one sub-problem could provide the possibility to eliminate many possible solutions from consideration. In other words, once a current candidate solution is obtained, we judge other nodes for further consideration. If the best potential solution of one particular node is no better than the current candidate solution, then there is no need to branch on this node. Otherwise, branching is needed to yield a new bound. In this section, we use one concrete example to illustrate how this method works.

Based on the theory we built in previous sections, here we propose an algorithm for finding an optimal solution of problem (1).

**Step 1:** Find the maximum solution of system (1). Compute $\bar{X} = \min_{i \in I} \{\hat{x}_i\}$ according to relation (2).

**Step 2:** Check feasibility. If $A\varphi\bar{X} = b$, continue. Otherwise, stop! $S(A, b) = \emptyset$, and problem (1) has no feasible solution.

**Step 3:** Compute index sets. Compute define $\bar{J}_i = \{j \in J_i : T_L(a_{ij}, \bar{X}_j) = b_i\}, \forall j \in J$.

**Step 4:** Arrange cost vector. Define $c_j^+$ and $c_j^-$ for $j = 1, 2, \ldots, n$ and define sub-problems $Z_1$ and $Z_2$.

**Step 5:** Use the branch-and-bound concept based on $\{\underline{X}(e) : e \in \bar{E}\}$ to solve $Z_1$.

**Step 6:** Generate an optimal solution of (1) via formula (5).

Consider problem (1) as follows:

$$\min 3x_1 + 4x_2 + x_3 + x_4 - x_5 + 5x_6$$

$$\begin{bmatrix} 0.5 & 0.85 & 0.9 & 0.3 & 0.85 & 0.4 \\ 0.2 & 0.2 & 0.1 & 0.95 & 0.1 & 0.8 \\ 0.8 & 0.8 & 0.4 & 0.1 & 0.1 & 0.1 \\ 0.1 & 0.1 & 0.1 & 0.1 & 0.1 & 0 \end{bmatrix} \varphi x = \begin{bmatrix} 0.85 \\ 0.6 \\ 0.5 \\ 0.1 \end{bmatrix}$$

**Step 1:** We find $\bar{X} = [0.7, 0.7, 0.95, 0.65, 1, 0.8]$.

**Step 2:** Since $A\varphi\bar{X} = b$, we know $S(A, b) \neq \emptyset$.

**Step 3:** We compute $J_1 = \{2,3,5\}$, $J_2 = \{4,6\}$, $J_3 = \{1,2\}$ and $J_4 = \{1,2,3,4,5\}$. So, the cardinality of $E$ and $\{\underline{X}(e) : e \in E\}$ is equal to 60. However, in this example we have also $\bar{J}_1 = \{3,5\}$, $\bar{J}_2 = \{4,6\}$, $\bar{J}_3 = \{1,2\}$ and $\bar{J}_4 = \{5\}$. Therefore, the cardinality of $\bar{E}$ and $\{\underline{X}(e) : e \in \bar{E}\}$ is equal to 8.

**Step 4:** The main problem is converted to the following two sub-problems:

$$\min Z_1 = 3x_1 + 4x_2 + x_3 + x_4 + 5x_6$$

$$\begin{bmatrix} 0.5 & 0.85 & 0.9 & 0.3 & 0.85 & 0.4 \\ 0.2 & 0.2 & 0.1 & 0.95 & 0.1 & 0.8 \\ 0.8 & 0.8 & 0.4 & 0.1 & 0.1 & 0.1 \\ 0.1 & 0.1 & 0.1 & 0.1 & 0.1 & 0 \end{bmatrix} \varphi x = \begin{bmatrix} 0.85 \\ 0.6 \\ 0.5 \\ 0.1 \end{bmatrix}$$

$$\min Z_2 = -x_5$$

$$\begin{bmatrix} 0.5 & 0.85 & 0.9 & 0.3 & 0.85 & 0.4 \\ 0.2 & 0.2 & 0.1 & 0.95 & 0.1 & 0.8 \\ 0.8 & 0.8 & 0.4 & 0.1 & 0.1 & 0.1 \\ 0.1 & 0.1 & 0.1 & 0.1 & 0.1 & 0 \end{bmatrix} \varphi x = \begin{bmatrix} 0.85 \\ 0.6 \\ 0.5 \\ 0.1 \end{bmatrix}$$

**Step 5:** Our branch-and-bound method begins by choosing j from $\bar{J}_1 = \{3,5\}$. Hence, each feasible solution must satisfy either $e(1) = 3$ or $e(1) = 5$. This yields two branches denoted by nodes 1 and 2 in Fig. 1. We then obtain a lower bound on the Z-value associated with each node. For example, when $e(1) = 3$, we know $c_3(b_1 + 1 - a_{13}) = 1(0.95) = 0.95$, therefore, any choice with $e(1) = 3$ results in $Z \geq 0.95$. So we write $Z \geq 0.95$ for node 1 of Fig. 1. Similarly, any choice with $e(1) = 5$ results in $Z \geq 0$. Since at this moment we have no reason to exclude any one of nodes 1 and 2 from consideration, both nodes have to be further investigated. By using the jump-tracking technique, we branch on the node with a lower bound on Z. In this case, we choose node 2. Since $\bar{J}_2 = \{4,6\}$, any choice associated with node 2 must satisfy either $e(2) = 4$ or $e(4) = 6$. Branching on node 2 yields nodes 3 and 4 in Fig. 1. For each new node, we need to evaluate a lower bound for the objective value Z. For example, at node 3, we compute $Z \geq c_5(b_1 + 1 - a_{15}) + c_4(b_2 + 1 - a_{24}) = 0.65$. Hence, any choice associated with node 3 must have $Z \geq 0.65$. Similar reasoning shows that $Z \geq 4$ for node 4. Again, we do not have any reason to exclude any of nodes 3 and 4 from consideration, so we need to branch on one node. The jump-tracking technique directs us to branch node 3.

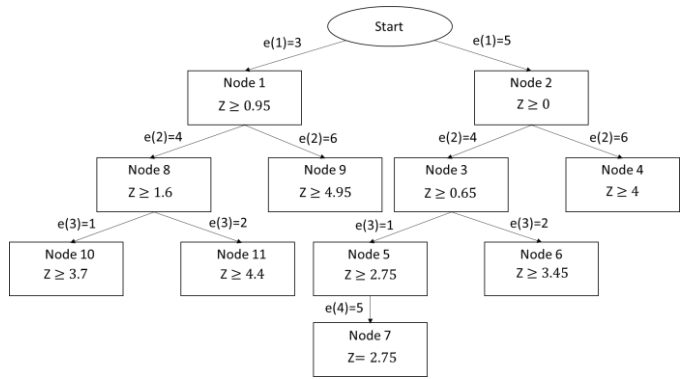

Figure 1. Modified branch and bound method.

Any choice associated with node 3 must satisfy either $e(3) = 1$ or $e(3) = 2$, since $\bar{J}_3 = \{1,2\}$. This yields nodes 5 and 6 in Fig. 1. By the same reasoning as before, any choice associated with node 5 must have $Z \geq 2.75$ and any choice associated with node 6 must have $Z \geq 3.45$. Of course, we are interested in node 5. To branch further on node 5, any choice associated with node 5 must have $e(4) = 5$, since $\bar{J}_4 = \{5\}$. This yields node 7 in Fig. 1. Note that node 7 corresponds to the sequence $e = [5,4,1,5]$. This sequence leads to a value of $Z = 2.75$. Therefore, node 7 is a feasible sequence which may be viewed as a candidate solution with $Z = 2.75$. Because the Z value for nodes 4 and 6 cannot be lower than 2.75, these nodes can be eliminated from further consideration. However, node 1 cannot be eliminated yet, because it is still possible for node 1 to yield a sequence having $Z < 2.75$. Hence, we now branch on node 1. Since $\bar{J}_2 = \{4,6\}$, any sequence associated with node 1 must have either $e(2) = 4$ or $e(4) = 6$. Correspondingly, we create nodes 8 and 9, respectively. For node 9, we calculate $Z \geq 4.95$. Since any sequence associated with node 9 must have $Z \geq 4.95 > 2.75$, this node is eliminated from consideration. For node 8, we have $Z \geq 1.6$. So it cannot be eliminated yet. We branch on node 8. Since $\bar{J}_3 = \{1,2\}$, any sequence associated with node 8 must have either $e(3) = 1$ or $e(3) = 2$. This creates nodes 10 and 11 in Fig. 1. For node 10, since $Z \geq 3.7$, we eliminate it. Similarly, for node 11 with $Z \geq 4.4$, it can be eliminated too. Now, with the exception of node 7, every node in Fig. 1 has been eliminated from consideration. Node 7 yields the sequence $e = [5,4,1,5]$ that generates solution $X(e) = [0.7, 0, 0, 0.65, 1, 0]$. Thus, we should choose the sequence $e = [5,4,1,5]$ and solution $X(e) = [0.7, 0, 0, 0.65, 1, 0]$ with $Z = 2.75$ as an optimal solution for $Z_1$.

**Step 6:** According to formula (5), we can find the optimal value of the original problem as $x^* = [0.7, 0, 0, 0.65, 1, 0]$ with optimal objective value $Z^* = 2.75$.

## V. NUMERICAL EXAMPLES

In this section, we present the experimental results for evaluating the performance of our algorithm. For this purpose, we apply our algorithm to the following eight test problems. The test problems have been randomly generated in different sizes.

**Test Problem 1:**

$c = [-9.3708, 6.7080, 6.7142, -9.002, 0.9177, 8.8633]$

$$\begin{bmatrix} 0.3214 & 0.7991 & 0.4896 & 0.2989 & 0.8159 & 0.8991 \\ 0.8064 & 0.0495 & 0.9728 & 0.2561 & 0.0983 & 0.8999 \\ 0.6013 & 0.2831 & 0.7484 & 0.8865 & 0.8595 & 0.5241 \\ 0.7896 & 0.6534 & 0.5678 & 0.4468 & 0.0276 & 0.1201 \end{bmatrix} \varphi x = \begin{bmatrix} 0.3632 \\ 0.3620 \\ 0.1490 \\ 0.2641 \end{bmatrix}$$

**Test Problem 2:**

$c = [1.337\ 3.6078\ -2.5724\ -8.4354\ -0.8729\ -9.0431]$

$$\begin{bmatrix} 0.7382 & 0.9374 & 0.7589 & 0.1100 & 0.2611 & 0.1320 \\ 0.0380 & 0.5133 & 0.9933 & 0.5970 & 0.0948 & 0.4528 \\ 0.9542 & 0.2409 & 0.3567 & 0.4305 & 0.4509 & 0.6521 \\ 0.7423 & 0.2599 & 0.7528 & 0.7307 & 0.6400 & 0.8269 \end{bmatrix} \varphi x = \begin{bmatrix} 0.1641 \\ 0.2724 \\ 0.0316 \\ 0.3083 \end{bmatrix}$$

**Test Problem 3:**

$c = [-0.0655\ 6.1730\ 2.6573\ 3.7680\ 2.7914]$

$$\begin{bmatrix} 0.7293 & 0.1635 & 0.4752 & 0.1486 & 0.1663 \\ 0.8598 & 0.9060 & 0.8053 & 0.6581 & 0.1496 \\ 0.6269 & 0.0773 & 0.5307 & 0.6339 & 0.2027 \\ 0.1805 & 0.3385 & 0.2273 & 0.2293 & 0.9549 \\ 0.5733 & 0.5806 & 0.7094 & 0.18222 & 0.0159 \end{bmatrix} \varphi x = \begin{bmatrix} 0.0325 \\ 0.2973 \\ 0.2725 \\ 0.0987 \\ 0.1177 \end{bmatrix}$$

**Test Problem 4:**

$c = [-7.798\ -3.1243\ -4.7186\ 1.6036\ -1.2921]$

$$\begin{bmatrix} 0.4151 & 0.6588 & 0.5561 & 0.6565 & 0.0458 \\ 0.0542 & 0.1191 & 0.8343 & 0.2774 & 0.1013 \\ 0.2596 & 0.1292 & 0.1061 & 0.7091 & 0.4689 \\ 0.2441 & 0.0779 & 0.2565 & 0.1302 & 0.7106 \\ 0.7090 & 0.7352 & 0.0541 & 0.5239 & 0.6013 \end{bmatrix} \varphi x = \begin{bmatrix} 0.3601 \\ 0.4373 \\ 0.3676 \\ 0.3233 \\ 0.5417 \end{bmatrix}$$

**Test Problem 5:**

$c = [0.8080\ 7.0530\ -5.5483\ 2.7401\ 4.2496\ -3.3253]$

$$\begin{bmatrix} 0.0507 & 0.7482 & 0.5739 & 0.3970 & 0.4639 & 0.2466 \\ 0.4593 & 0.6247 & 0.8003 & 0.3888 & 0.7690 & 0.0983 \\ 0.8656 & 0.7280 & 0.4496 & 0.6059 & 0.3349 & 0.4098 \\ 0.7758 & 0.6701 & 0.2760 & 0.0185 & 0.1620 & 0.8231 \\ 0.3420 & 0.4851 & 0.7258 & 0.6908 & 0.5847 & 0.1639 \end{bmatrix} \varphi x = \begin{bmatrix} 0.0935 \\ 0.1484 \\ 0.4456 \\ 0.4797 \\ 0.1326 \end{bmatrix}$$

**Test Problem 6:**

$c = [7.6046\ 6.1145\ 1.9805\ -6.2986\ 8.6347\ 3.1955]$

$$\begin{bmatrix} 0.6940 & 0.9782 & 0.4444 & 0.1527 & 0.1523 & 0.692 \\ 0.8041 & 0.2452 & 0.1169 & 0.3113 & 0.5080 & 0.5420 \\ 0.3796 & 0.3927 & 0.9970 & 0.5338 & 0.7253 & 0.1315 \\ 0.2498 & 0.0264 & 0.6789 & 0.4315 & 0.7274 & 0.4455 \\ 0.3284 & 0.7754 & 0.4421 & 0.8913 & 0.8067 & 0.4915 \end{bmatrix} \varphi x = \begin{bmatrix} 0.0748 \\ 0.1226 \\ 0.1888 \\ 0.1415 \\ 0.2296 \end{bmatrix}$$

**Test Problem 7:**

$c = [1.6779\ -9.7153\ -1.8290\ -2.3325\ -7.1269\ 4.7970]$

$$\begin{bmatrix} 0.3247 & 0.2323 & 0.1009 & 0.8162 & 0.8712 & 0.0202 \\ 0.3829 & 0.9656 & 0.5807 & 0.5737 & 0.7469 & 0.1203 \\ 0.9890 & 0.1448 & 0.1044 & 0.2106 & 0.3966 & 0.6341 \\ 0.6230 & 0.4010 & 0.5958 & 0.4310 & 0.1664 & 0.4847 \\ 0.4039 & 0.5320 & 0.4781 & 0.4327 & 0.5418 & 0.9948 \end{bmatrix} \varphi x = \begin{bmatrix} 0.3249 \\ 0.1082 \\ 0.0229 \\ 0.0830 \\ 0.0250 \end{bmatrix}$$

**Test Problem 8:**

$c = [-6.1413\ -7.2350\ -0.3884\ -8.7017\ 7.3308\ -6.7394\ -1.9768]$

$$\begin{bmatrix} 0.7073 & 0.6762 & 0.1256 & 0.3776 & 0.3128 & 0.9969 & 0.9627 \\ 0.4065 & 0.5206 & 0.5880 & 0.9598 & 0.1945 & 0.4147 & 0.0247 \\ 0.7558 & 0.9522 & 0.2166 & 0.4235 & 0.1167 & 0.5126 & 0.0011 \\ 0.8589 & 0.7442 & 0.0344 & 0.1647 & 0.9460 & 0.1518 & 0.5316 \\ 0.8385 & 0.7224 & 0.9988 & 0.2616 & 0.4401 & 0.4538 & 0.7876 \end{bmatrix} \varphi x = \begin{bmatrix} 0.1741 \\ 0.0950 \\ 0.0017 \\ 0.3869 \\ 0.1432 \end{bmatrix}$$

Table 1 shows the cardinality of $E$ and $\{\underline{X}(e) : e \in E\}$, the cardinality of $\overline{E}$ and $\{\underline{X}(e) : e \in \overline{E}\}$, the number of nodes that the algorithm visits and the number of complete paths visited by the algorithm.

Table 1. The results found by our algorithm.

| Test problems | Cardinality of $E$ | Cardinality of $\overline{E}$ | Number of visited nodes | Number of visited paths |
|---|---|---|---|---|
| 1 | 288 | 4 | 6 | 1 |
| 2 | 480 | 3 | 6 | 1 |
| 3 | 1200 | 1 | 5 | 1 |
| 4 | 24 | 1 | 5 | 1 |
| 5 | 1800 | 2 | 7 | 1 |
| 6 | 4500 | 2 | 9 | 1 |
| 7 | 2592 | 2 | 9 | 2 |
| 8 | 6048 | 4 | 13 | 1 |

Moreover, for each test problem $i = 1, ..., 8$, the best path $e_i$, the optimal solution $x_i^*$ and the optimal value $z_i^*$ is listed as follows:

$e_1 = [2\ 3\ 4\ 1]$

$x_1^* = [0.4745\ 0.5641\ 0.3892\ 0.2625\ 0\ 0]$

$z_1^* = -0.41232$

$e_2 = [2\ 3\ 5\ 4]$

$x_2^* = [0\ \ 0.2267\ \ 0.2791\ \ 0.5776\ \ 0.5807\ \ 0.3795]$

$z_2^* = -8.7111$

$e_3 = [1\ 2\ 4\ 5\ 3]$

$x_3^* = [0.3032\ \ 0.3913\ \ 0.4082\ \ 0.6386\ \ 0.1437]$

$z_3^* = 6.2883$

$e_4 = [2\ 3\ 4\ 5\ 1]$

$x_4^* = [0.8326\ \ 0.7013\ \ 0.6030\ \ 0.6585\ \ 0.6127]$

$z_4^* = -11.2659$

$e_5 = [2\ 3\ 1\ 6\ 4]$

$x_5^* = [0.5800\ \ 0.3453\ \ 0.3480\ \ 0.4417\ \ 0\ \ 0.6565]$

$z_5^* = 0.0008$

$e_6 = [2\ 1\ 3\ 5\ 4]$

$x_6^* = [0.3181\ \ 0.0966\ \ 0.1917\ \ 0.3382\ \ 0.4141\ \ 0]$

$z_6^* = 4.8385$

$e_7 = [4\ 2\ 1\ 3\ 6]$

$x_7^* = [0.0338\ \ 0.1425\ \ 0.4871\ \ 0.5087\ \ 0.3612\ \ 0.0302]$

$z_7^* = -5.8352$

$e_8 = [6\ 4\ 1\ 5\ 3]$

$x_8^* = [0.2459\ \ 0.0495\ \ 0.1444\ \ 0.1352\ \ 0.4409\ \ 0.1772\ \ 0.2113]$

$z_8^* = -1.4818$

## VI. Conclusion

In this paper, we have studied a linear optimization problem subjected to a system of Lukasiewicz fuzzy relation equations and presented a procedure to find an optimal solution. Due to the non-convexity nature of its feasible domain, we tend to believe that there is no polynomial-time algorithm for this problem. The best we can do here is that, after analyzing the properties of its feasible domain, we find a solution set that include all the minimal solutions of the main feasible region, then apply a modified branch-and-bound method to find one solution. The question of how to generate the whole optimal solution set is yet to be investigated. From the analysis of Theorem 2, it is clearly seen that if all minimum solutions of a given system of fuzzy relation equations can be found, then an optimal solution of the optimization problem defined by (5) can be constructed. Therefore, solving this optimization problem is no harder than solving a system of fuzzy relation equations for all minimum solution. Although it is not known whether these two problems are essentially equivalent or not, the basic concepts introduced in Section 2 have been further developed for solving fuzzy relation equations. Extension to other types of objective functions is under investigation.